\documentclass[11pt,a4paper]{article}
\newcommand{\bm}[1]{\mbox{\boldmath $#1$}}
\usepackage{epsfig}
\usepackage {amsfonts}
\usepackage {amsmath}
\usepackage{latexsym}

\begin{document}
\title{Adaptive computation of the Symmetric Nonnegative Matrix
Factorization (NMF)}
\author {P. Favati \and G. Lotti \and O. Menchi \and F. Romani}

\date { }

\maketitle
\begin{abstract}
Nonnegative Matrix
Factorization (NMF), first proposed in 1994 for data
analysis, has received successively much attention in a great variety of
contexts such as data mining, text clustering, computer vision, bioinformatics, etc.
In this paper the case of a symmetric matrix is considered and the symmetric nonnegative
matrix factorization (SymNMF) is obtained by using a penalized nonsymmetric
minimization problem. Instead of letting the penalizing parameter increase according to an a priori fixed rule, as
suggested in literature,
 we propose a heuristic approach
based on an adaptive technique. Extensive experimentation shows that the proposed algorithm is effective.

\end{abstract}

\section{Introduction}\label{intro}
Dimensional reduction problems are
of fundamental relevance for data compression and
classification. An important problem of this kind is represented by the
nonnegative matrix factorization, which was first proposed in \cite{paa} for data
analysis and afterwards widely applied (see \cite{kimpark} for an extensive bibliography).

Let ${\bf R}_+^m$ be the $m$-dimensional space of vectors with nonnegative
components and $M$ a matrix of $n$ columns ${\bm m}_i \in {\bf R}_+^m$, for
$i=1,\ldots, n$. Given an integer $k<\min(m,n)$, the problem of finding two
low-rank matrices $W\in {\bf R}_+^{m\times k}$ (the {\it basis} matrix) and
$H\in {\bf R}_+^{n\times k}$ (the {\it coefficient} matrix) such that the
product $WH^T$ approximates  $M$, is known as {\it Nonnegative Matrix
Factorization} (NMF). In this way the $n$ objects ${\bm m}_i$ result represented by linear combinations with nonnegative coefficients of few nonnegative basis vectors.

A specific formulation of the problem requires that a metric is assigned to
measure the distance between $M$ and $WH^T$. The
nonnegativity of the involved items would suggest to minimize a divergence, like the likelihood {\it Kullback-Leibler} divergence, but some computational
difficulties and the slow convergence rate of common iterative
procedures used to tackle the problem, suggest the more flexible
metric of the F-norm (Frobenius norm). With this norm the general NMF problem takes the form
 \begin{equation}\label{pro}
\min_{W,\,H\ge O}\,
 \Phi(M,W,H),\ \hbox{where}\
 \Phi(M,W,H)=\textstyle\frac12\;\|M-W\,H^T\|_F^2.
 \end{equation}

The best factorization of a matrix in terms of the F-norm is
achieved by the Singular Value Decomposition $M=U{\mit \Sigma} V^T$.
Then the best $k$-rank approximation of $M$ is $ U_k {\mit \Sigma}_k
V_k^T$, where $U_k$ and $ V_k$ are the truncated submatrices of $U$
and $V$ to $k$ columns and ${\mit \Sigma}_k$ is the $k\times k$
leading submatrix of ${\mit \Sigma}$. Unfortunately, nothing can be
said about the sign of the entries of $ U_k$ and $ V_k$, which could
be negative. It follows that nonnegativity must be imposed as a
constraint. Other constraints could appear in (\ref{pro}) in order to satisfy additional requirements \cite{ppp}.
 For example, adding the  terms $\rho_1 \|W\|+\rho_2 \|H\|$ to the function $\Phi$ to be minimized, where  $\|\cdot \|$  is a  suitable norm and  $\rho_1$ and $\rho_1$ are positive parameters, would give a regularized solution and possibly control the sparsity of the factors $W$ and $H$.

In this paper we consider the additional requirement of symmetry:
problem (\ref{pro}) is then replaced by the symmetric NMF ({\it SymNMF}) problem
\begin{equation}\label{sympro}
  \min_{W\ge O}\,
 \Psi(A,W),\ \hbox{where}\
 \Psi(A,W)=\textstyle\frac12\;\|A-W\,W^T\|_F^2,
 \end{equation}
where $A$ is a symmetric $n\times n$ matrix of nonnegative entries and $W\in {\bf R}_+^{n\times k}$. Note that the required approximation $WW^T$ is positive semidefinite and can be very poor if $A$ does not have enough nonnegative eigenvalues.

Problem (\ref{sympro}) has a fourth-order nonconvex  objective function and
optimization algorithms guarantee only the stationarity of the limit points, so
one only looks for a local minimum. Standard gradient algorithms lead to
stationary solutions, but suffer from slow convergence. Newton-like algorithms,
which have a better rate of convergence, are computationally expensive. In
\cite{kuang1} a nonsymmetric formulation of (\ref{sympro}) is suggested by
considering the following penalized problem
 \begin{equation}\label{nonsympro}
 \min\limits_{W,\,H\ge O}\,
 \Phi_\alpha(A,W,H), \ \hbox{where}\
 \Phi_\alpha(A,W,H)=\textstyle\frac12\;\Big(
 \|A-W H^T\|_F^2+\alpha  \|W-H \|_F^2\Big),
 \end{equation}
$\alpha$ being a positive parameter which acts on the violation of the
symmetry. Choosing  $\alpha$  aligned with the magnitude of $A$, makes the penalized problem invariant from the
scale of matrix $A$. In this paper we propose an algorithm to approximate the solution of (\ref{sympro}) by solving
 iteratively problem (\ref{nonsympro}) and dealing adaptively with the penalizing parameter $\alpha$.

To this aim in Section \ref{anls} we recall the ANLS framework, a standard approach for tackling a general
 (i.e. not symmetric) NMF problem of form (\ref{pro}) by addressing alternatively two convex subproblems,
together with the two methods ({\tt BPP} \cite{kimpark2} and {\tt GCD} \cite{dhil}) which will be used
to solve each subproblem.  Section \ref{snmf} deals with the heuristic for the choice of the parameter
$\alpha$  to solve (\ref{nonsympro}). In Section \ref{expe} the results of an extensive experimentation are presented to validate the proposed
adaptive strategy and to compare the performance of the two chosen methods when applied in our context.

\section{The ANLS framework for the general NMF}\label{anls}
Problem (\ref{pro}) is nonconvex and finding its global minimum is
NP-hard. Most nonconvex optimization algorithms guarantee only the
stationarity of the limit points, so one looks for a local
minimum. There is a further source of nonunicity, since
$WH^T=W'{H'}^T$ with $W'=WS$, $H'=HS^{-T}$, where $S\in {\bf R}_+^{k
\times k}$ is a nonsingular scaling matrix. This can be fixed by
choosing for example $S$ in such a way to normalize the columns of
$W$ to unit 2-norm.

The {\it alternating nonnegative least squares} (ANLS)
method, which belongs to the {\it block coordinate descent} (BCD) framework of
nonlinear optimization \cite{kimpark3}, solves iteratively problem (\ref{pro}). First,
one of the factors, say $W$, is initialized to $W^{(0)}$ with
nonnegative entries and the matrix $H^{(1)}\in {\bf R}_+^{n\times k}$
realizing the minimum of
$\Phi(M,W^{(0)},H)$ on ${H \ge 0}$ is computed. Then a new matrix $W^{(1)}\in {\bf
R}_+^{m\times k}$ realizing the minimum of $\Phi(M,W,H^{(1)})$ on ${W \ge 0}$ is computed, and so on, updating $W$ and $H$ alternatively.
In practice the following {\it inner-outer} scheme is
applied
\begin{align}
 H^{(\nu)}=\mathop{\rm argmin}\limits_{H \ge 0}\,
 \Phi(M,W^{(\nu-1)},H),\label{proh}
 \\
 W^{(\nu)}=\mathop{\rm argmin}\limits_{W \ge 0}\,
 \Phi(M^T,H^{(\nu)},W),\label{prow}
 \end{align}
for $\nu=1,2,\ldots$, where each subproblem is solved by applying a chosen inner method.
At the $\nu$th outer iteration a suitable stopping condition should check whether
a local minimum of the object function $\Phi(M,W,H)$ of (\ref{pro})
has been sufficiently well approximated, for example by monitoring the {\it error}, i.e.
the distance of $W^{(\nu)}H^{(\nu)T}$ from $M$
\[
 e^{(\nu)}=\|M-W^{(\nu)}H^{(\nu)T}\|_F^2.
 \]
The choice of the initial matrix $W^{(0)}$ may be critical, due to
the fact that only a local minimum is expected, which obviously
depends on this choice  and, typically, the algorithm is run several
times with different initial matrices.

Although the original problem (\ref{pro}) is nonconvex, subproblems
(\ref{proh}) and (\ref{prow}) are convex and nearly any procedure for constrained quadratic
optimization can be chosen as inner method (for example an
Active-Set-like method \cite{bjo, kimpark2, law}).  The requirement that the inner problems
are exactly solved at each outer step is necessary for convergence
\cite{grippo} but makes the overall algorithm rather slow at large
dimensions. Faster
approaches have been devised by computing iteratively approximate
solutions with inexact methods like modified gradient descent methods or projected Newton-type methods \cite{kimpark3}. In this paper we take into consideration,
as inner methods,  an Active-Set-like method with block principal pivoting (the {\tt BPP} method, coded as Algorithm 2 in \cite{kimpark2}) and a coordinate
descent method, called Greedy Coordinate Descent
({\tt GCD}) in \cite{dhil}). Their main difference lies in the termination: exact for {\tt BPP} and approximated for {\tt GCD}.

When the ANLS method is applied, at each outer step,
say the $\nu$th outer step, the inner method computes the solution of two problems of the form
\begin{equation}\label{inme}
 \min_{X \ge 0}\,\Phi(B,C,X)=\min_{X \ge 0}\,
 \textstyle\frac12\;\|B-C\,X^T\|_F^2,
 \end{equation}
where $B=M$, $C=W^{(\nu-1)}$, $X=H$ for problem (\ref{proh})
and $B=M^T$, $C=H^{(\nu)}$, $X=W$ for problem (\ref{prow}). We assume matrix $C$ to have full rank.

Let $r\times s$ be the dimensions of $B$  ($r=m$, $s=n$ in the first case and $r=n$, $s=m$ in the second case). Denoting by ${\bm b}\in {\bf R}_+^r$ and ${\bm x}\in {\bf R}^k$ the $h$th columns of $B$ and $X^T$ respectively, for $h=1,\ldots,s$, problem (\ref{inme}) can be decomposed into $s$ independent
least squares nonnegatively constrained problems
 \begin{equation}\label{sis2}
 \min_{{\bm x} \ge 0}\varphi({\bm x}),\quad\hbox{where}\quad \varphi({\bm x})=\textstyle\frac12\;\|
 \,{\bm b}-C{\bm x}\,\|_2^2.
 \end{equation}
 The gradient of the objective function $\varphi({\bm x})$ is
${\bm g}({\bm x})=C^T(C{\bm x}-{\bm b})$.

The $s$ problems (\ref{sis2}) are solved in sequence, using either {\tt BPP} or {\tt GCD}.
Before proceeding, we give a brief description of the two considered methods. The corresponding
codes can be found in the cited papers.

\subsection{The {\tt BPP} method}\label{inner1}
{\tt BPP} method derives from the classical active set method for linearly constrained optimization. For a point ${\bm x}\in
{\bf R}^k$ consider the active and passive index sets
at ${\bm x}$
 \[
   {\cal A}({\bm x})=\{1\le i\le k,\ \hbox{such that}\ x_i=0\},
   \quad
   {\cal P}({\bm x})={\cal K}\,-\,{\cal A}({\bm x}),\]
where ${\cal K}=\{1,\ldots,k\}$ is the complete index set.
Let $C_{\cal A}$ and $C_{\cal P}$ be the restrictions of the matrix $C$ to ${\cal A}({\bm x})$ and
${\cal P}({\bm x})$ respectively.
Since $C_{\cal P}$ has full column rank, the solution of the unconstrained least
squares problem
\begin{equation}\label{unc0}
   {\bm x}_{\cal P}=
 \mathop{\rm argmin}\limits_{{\bm z}}
 \textstyle\frac12\;\|\,{\bm b}-C_{\cal P}\,{\bm z}\,\|_2^2
\end{equation}
is given by the solution
of the system
 \begin{equation}\label{unc1}
 C_{\cal P}^TC_{\cal P}\,{\bm z}=C_{\cal P}^T{\bm b},
 \end{equation}
which has size less than or equal to $k$.
If the size is not too large, the system is solved by applying the Cholesky factorization (otherwise, one can resort to the conjugate gradient).
Let ${\bm x}^*$ be the vector which coincides with ${\bm x}_{\cal P}$ on ${\cal P}({\bm x})$ and has zero components on ${\cal A}({\bm x})$. Denote by
 \begin{equation}\label{unc2}
 {\bm g}_{\cal A}({\bm x}^*)=C_{\cal A}^T(C_{\cal P}{\bm x}_{\cal P}-{\bm b})
  \end{equation}
the gradient restricted to ${\cal A}({\bm x})$.
According to Karush-Kuhn-Tucker (KKT) optimality conditions, the
vector ${\bm x}^*$ is a solution of (\ref{sis2}) if and
only if ${\bm x}_{\cal P}\ge {\bm 0}$ and ${\bm g}_{\cal A}({\bm x}^*) \ge {\bm 0}$.

If the active and passive index sets of ${\bm x}^*$ were known in advance, problem
(\ref{sis2}) could be solved by simply solving (\ref{unc1}). Since
the two index sets are initially unknown,
a sequence of unconstrained subproblems is solved with the two index sets
 in turn predicted
and exchanged.  The computation starts with index sets associated to an initial point supplied by the outer iteration and  goes on until all the constraints
become passive or the gradient has nonnegative components
corresponding to all the active constraints, indicating that the
objective function cannot be reduced any more.

In the classical Active-Set method \cite{law}, only an index moves from ${\cal A}({\bm x})$ to ${\cal
P}({\bm x})$ at a time. This makes the number of iterations to grow
considerably with the size of the problem. An overcome to this
drawback consists in exchanging more indices between ${\cal A}({\bm
x})$ and ${\cal P}({\bm x})$ at each iteration, as suggested in
\cite{kimpark2}. The number of iterations results reduced, but the
generated vectors ${\bm x}$ may fail to maintain nonnegativity and
the monotonic decrease of the objective function is not guaranteed.
A finite termination is achieved by a backup rule which implements
the standard one index exchange when necessary.

When the procedure described above for a single column ${\bm b}$ of $B$ is applied to all the columns of $B$, the following  improvement,  proposed in \cite{kimpark2}, reduces the computational cost. Since  each problem (\ref{sis2}) shares the same matrix $C$, and the main cost depends on solving system (\ref{unc1}) and on computing vector (\ref{unc2}) with matrices $C_{\cal P}^TC_{\cal P}$, $C_{\cal A}^TC_{\cal P}$  and vectors $C_{\cal P}^T{\bm b}$, $C_{\cal A}^T{\bm b}$, it is
suggested to extract these matrices and vectors from the complete matrices $C^TC$ and $C^TB$ computed once at the beginning.
Another improvement consists in grouping the right-hand side vectors which share the same index set ${\cal P}$ in order to avoid
repeated computation of the Cholesky factorization in solving systems (\ref{unc1}).

\subsection{The {\tt GCD} method}\label{inner2}

{\tt GCD} derives from {\tt FastHals} \cite{cic}, an iterative method which performs a cyclic coordinate descent scheme. {\tt GCD}, instead, at each step selectively replaces the element whose update leads to the largest decrease of the objective function.

In  \cite{dhil} {\tt GCD} works on the entire matrix $X$, but in practice the method is applied to solve in sequence problems of
type (\ref{sis2}). For each problem  (\ref{sis2}), starting from a ${\bm x}^{(0)}\in {\bf R}_+^k $ chosen according to the outer
iteration, {\tt GCD} computes a sequence ${\bm x}^{(j)}$,  $j=1,2,\ldots$,  until suitably stopped. A global stopping condition
based on the entire matrix $B$, suggested in \cite{dhil}, is described at the end of the paragraph.

At the $j$th iteration the vector
${\bm x}^{(j)}$ is obtained by applying a single coordinate correction according
to the rule
 \[
 {\bm x}^{(j)}= {\bm x}^{(j-1)}+\widehat\lambda\,{\bm e}_i,
 \]
where $i$ is an index to be  selected in $\{1,\ldots,k\}$, ${\bm e}_i$ is the $i$th canonical $k$-vector
and the scalar $\widehat\lambda$ is determined by imposing that
$\varphi({\bm x}^{(j)})$, as a function of $\lambda$, is the minimum on the set
$
  S=\{\lambda\ \hbox{such that}\ x_i^{(j-1)}+\lambda\ge 0\}.
$
This scalar is
 \[
 \widehat\lambda=
    -\,\frac{g_i^{(j-1)}}{q_{i,i}}\quad
    \hbox{if}\quad \displaystyle\frac{g_i^{(j-1)}}{q_{i,i}}\le
    x_i^{(j-1)}\quad\hbox{and}\quad
  \widehat\lambda=-\,x_i^{(j-1)}\quad\hbox{otherwise},
 \]
where $ Q=C^TC$ is the Hessian of $\varphi$
and ${\bm g}^{(j-1)}={\bm g}\big({\bm x}^{(j-1)}\big)$.
In correspondence, the objective function is decreased by
\[
 d_i^{(j)}=\varphi\big({\bm x}^{(j-1)}\big)-\varphi\big({\bm x}^{(j)}\big)=-\,g_i^{(j-1)}\,\widehat\lambda
 -\frac{1}{2}\,q_{i,i}\,\widehat\lambda^2.\]
A natural choice for index $i$ is
the one that maximizes $d_i^{(j)}$ varying $i$. As a consequence,  $x_i^{(j-1)}$ is updated by adding $\widehat\lambda$
and the elements of the gradient become
\[
g_t^{(j)}=g_t^{(j-1)}+\widehat\lambda \,q_{t,i},\quad\hbox{for}\quad
t=1,\ldots,k.
\]
Then a new iteration begins, where a new index $i$ is detected,
and so on, until a stopping condition is met. In \cite{dhil} the following condition is suggested
\begin{equation}\label{mu}
 \max_i d_i^{(j)}<\eta\, \mu,
\end{equation}
where the quantity $\mu$ is the largest possible reduction of all  the
objective functions $\varphi$ of problems (\ref{sis2}) varying ${\bm b}$, that can be expected when a
single element is modified at
the first iteration and $\eta$ is a preassigned tolerance. Of course, the value of $\eta$ influences the
convergence of the outer method,
hence the overall computational cost. In \cite{dhil} $\eta=10^{-3}$ is suggested. We will examine this
question in Section \ref{expe}.

\section{The SymNMF problem}\label{snmf}

We turn now to the SymNMF problem (\ref{sympro}).
As anticipated, its solution is computed through the nonsymmetric formulation (\ref{nonsympro}), applying
ANLS as the outer algorithm, i.e.  by alternating the solution
of the two subproblems,

\begin{equation}\label{proh1}
 H^{(\nu)}=\mathop{\rm argmin}\limits_{H \ge 0}\;
 \Phi\bigg(\left[\begin{array}{c} A\\ \sqrt \alpha\;W^{(\nu-1)T}\end{array}\right],
 \left[\begin{array}{c} W^{(\nu-1)}\\ \sqrt \alpha\;I_k\end{array}\right],H
 \bigg),
 \end{equation}
 \begin{equation}\label{prow1}
 W^{(\nu)}=\mathop{\rm argmin}\limits_{W \ge 0}\;
 \Phi\bigg(\left[\begin{array}{c} A\\ \sqrt \alpha\;H^{(\nu)T}\end{array}\right],
 \left[\begin{array}{c} H^{(\nu)}\\ \sqrt \alpha\;I_k\end{array}\right],W
 \bigg).
 \end{equation}

The corresponding function ${\tt Sym\_ANLS}$ is shown in Figure \ref{f2}, where by
${\tt Inner\_solve}\,\big(B,C,X_0\big)$ we denote the function used to solve (\ref{inme}) employing a  method which
starts with initial iterate $X_0\ge 0$.
For problems (\ref{proh1}) and (\ref{prow1}) both {\tt BPP} and {\tt GCD}, used as inner methods, can be implemented
without explicitly forming the four block matrices.

At the $\nu$th outer iteration the stopping condition checks whether a local minimum of the object function $\Phi_\alpha(A,W,H)$
of (\ref{nonsympro}) has been sufficiently well approximated  by monitoring
$\epsilon_S^{(\nu)}$  and $\delta^{(\nu)}$, where
 \begin{equation}\label{sto1}
 \epsilon_S^{(\nu)}=\frac{\|A-W^{(\nu)}W^{(\nu)T}\|_F}{\|A\|_F}
 \end{equation}
measures the objective function of problem (\ref{sympro}) and
\begin{equation}\label{sto2}
 \delta^{(\nu)}=\frac{\|W^{(\nu)}-H^{(\nu)}\|_F}{\min(\|W^{(\nu)}\|_F,\|H^{(\nu)}\|_F)}
 \end{equation}
measures the {\it degree of symmetry}.

The starting points $W^{(0)}$ and $H^{(0)}$ are required by the first call of the inner method. Both {\tt BPP} and {\tt GCD}
 can start with $H^{(0)}=O$, because the gradient of the objective function
 in (\ref{proh1}) evaluated in the starting point is
\[G\big(H^{(0)}\big)=-W^{(0)T}(A+\alpha I_n)
\le O.\]
As $W^{(0)}$ we suggest a matrix of the form $R \sqrt{\|A\|_F}/\|R\|_F$ where $R$ is a random matrix with entries in $[0,1]$.
Moreover, the function ${\tt Sym\_ANLS}$ needs a procedure for updating the value of $\alpha$.

Let $(W_\alpha,H_\alpha)$ be the solution of (\ref{nonsympro}). The value of
$\alpha$  influences the symmetry of the solution: the largest $\alpha$, the
smallest $\|W_\alpha-H_\alpha\|$, but a too large $\alpha$, with respect to the magnitude of $A$, could lead to a
poor solution. If
$W_\alpha = H_\alpha$ for some  $\alpha$, then  $W_\alpha$ is also a solution of
(\ref{sympro}) and we call it a {\it symmetric} solution of (\ref{nonsympro}).
We call {\it quasi-symmetric} solution of (\ref{nonsympro}) a solution with
$W_\alpha\sim H_\alpha$ and ${\|A-W_\alpha H_\alpha^T\|^2}$ dominating over
$\alpha \|W_\alpha- H_\alpha\|^2$. In this case $W_\alpha$ is assumed as a good
approximation of the solution of (\ref{sympro}).

 If a quasi-symmetric solution  exists, it is possible that
the convergence to it is achieved even for a small $\alpha$,
provided that the starting point $W^{(0)}$ is sufficiently close to the quasi-symmetric solution. On
the other hand, a too large value of $\alpha$ should be avoided because
it tends to move the solution of (\ref{nonsympro}) away from a minimum point of (\ref{sympro}).

The sequence of penalizing  parameters is constructed by setting
\begin{equation}\label{ada1}
\alpha^{(\nu)}=\beta^{(\nu)}\, \max A, \quad\hbox{with} \quad \beta^{(0)}=1, \end{equation}
where  $\nu$ is the step index of the outer iteration.  The starting value $\alpha^{(0)}=\max A$ is
tuned according to the scale of $A$.

In  \cite{kuang1} the parameter $\beta$ is modified according to a geometric progression
of the form $\beta^{(\nu)}=\zeta^\nu$ where the fixed ratio $\zeta=1.01$ is suggested.
In the  experimentation the strategy proposed in \cite{kuang1} has been tested also  with different values  of $\zeta$.

Instead of considering an increasing sequence $\beta^{(\nu)}$, we suggest to modify the parameter $\beta^{(\nu)}$
in (\ref{ada1}) using an adaptive strategy,
called {\tt ADA}, which takes into account  the following quantities
\begin{equation}\label{sto3}
 \epsilon_N^{(\nu)}=\frac{\|A-W^{(\nu)}H^{(\nu)T}\|_F}{\|A\|_F}
 \end{equation}
which measures the first component of the objective function of problem
(\ref{nonsympro}), and the ratio
\begin{equation}\label{sto4}
  \rho^{(\nu)}=\frac {\epsilon_S^{(\nu)}}{\epsilon_N^{(\nu)}},
  \end{equation}
  where $\epsilon_S^{(\nu)}$ is defined in (\ref{sto1}).
One might think that $\epsilon_S^{(\nu)}\ge \epsilon_N^{(\nu)}$, but
this is not always true. In fact, if the stationary point to which
the outer method converges is a symmetric solution of problem  (\ref{nonsympro}), substituting
$H^{(\nu)}$ with $W^{(\nu)}$ can be seen as a sort of extrapolation,
that may even decrease the error.
\begin{figure} [ht]
 \hrule \vspace{0.1cm}
 \noindent{\tt function} ${\tt Sym\_ANLS}\,\big(W,H,\alpha,\nu_{max}\big)$\\
 {\sl computes
 recursively the solution of (\ref{sympro}) by solving (\ref{nonsympro}) given
 initial $W$, $H$, $\beta$ and the number of allowed
 iterations $\nu_{max}$.}\\[-0.2cm]
  \hrule
\begin{tabbing}
 \hspace*{0.5cm}\=$W^{(0)}=W$; \ $H^{(0)}=H$; \ $\beta^{(0)}=1$;
 \ $\nu=0$; \ $cond={\tt True}$;\\
 \>  {\tt whi}\={le}\ $cond$\\
 \> \> $\alpha^{(\nu)}=\beta^{(\nu)}\, \max A$;\\
 \> \>  $\nu=\nu+1$;\\
 \> \>  $\xi=\sqrt{\alpha^{(\nu-1)}}$;\\
 \> \>  $H^{(\nu)} ={\tt Inner\_solve}\,\bigg(\left[\begin{array}{c} A\\
  \xi \,W^{(\nu-1)T}\end{array}\right]
  ,\left[\begin{array}{c} W^{(\nu-1)}\\
  \xi\,I_k\end{array}\right],H^{(\nu-1)}\bigg);$\\
 \> \>  $W^{(\nu)} ={\tt Inner\_solve}\,\bigg(
  \left[\begin{array}{c} A\\ \xi\,H^{(\nu)T}\end{array}\right],
  \left[\begin{array}{c} H^{(\nu)}\\ \xi\,I_k\end{array}\right],
  W^{(\nu-1)}\bigg);$\\
  \> \>  compute $\epsilon_S^{(\nu)}$ and $\delta^{(\nu)}$,
  according to (\ref{sto1}) and (\ref{sto2});\\[0.1cm]
  \> \> $stop=\big|\epsilon_S^{(\nu)}- \epsilon_S^{(\nu-1)}
  \big|\le \tau_1\, \epsilon_S^{(\nu)}\ {\tt and}\ \delta^{(\nu)}\le\tau_2;$\\[0.1cm]
  \> \> $cond={\tt not}\ stop\ \ {\tt and}\ \ \nu<\nu_{max}$;\\
   \> \> $\beta^{(\nu)}={\tt Update}\,\big(\beta^{(\nu-1)}\big)$;\\[0.1cm]
 \> {\tt end while};\\[0.1cm]
 \> {\tt return} $W^{(\nu)}$;\\[-0.7cm]
\end{tabbing}
   \hrule
\caption{\label{f2} Algorithm to solve problem (\ref{sympro}). The function {\tt
Inner\_solve} solves problems (\ref{proh1}) and (\ref{prow1}). In all the experiments we have assumed
$\tau_1=10^{-3}$ and $\tau_2=0.1$.}
\end{figure}

When $\epsilon_S^{(\nu)}$ is  smaller than
$\epsilon_N^{(\nu)}$,
the value of $\beta^{(\nu)}$ can be safely reduced without risking an increase
of the distance from the symmetric solution. Otherwise the value of $\beta$ is increased depending on the value of $\rho^{(\nu)}$.
More precisely, when $\rho^{(\nu)}<1$ the decreasing rate of $\beta$ is tuned by the value of $\rho^{(\nu)}$, the degree of symmetry  $\delta^{(\nu)}$ and the magnitude of $\beta^{(\nu)}$.
Since in this case the outer iteration is well
directed towards a quasi-symmetric solution, the penalty condition can be relaxed without any risk.
When  $\rho^{(\nu)}>1$, the
value of $\beta$ is updated by means of multiplication by $\rho^{(\nu)2}$, paying
attention to avoid a too large increase. The adaptive strategy  is implemented  by function
{\tt ADA}, whose code is given in Figure \ref{f3}.

In Figure \ref{f2}, function {\tt Update} denotes the function used to update $\beta$. When the geometrical updating is chosen,
\[\beta ^{(\nu)}=\zeta \,\beta^{(\nu-1)}, \quad \hbox {for a fixed } \zeta.\]
When the adaptive updating is chosen,
\[\beta^{(\nu)}={\tt ADA}\,\big(W^{(\nu)},H^{(\nu)}, \delta^{(\nu)}, \epsilon_S^{(\nu)},\beta^{(\nu-1)}\big).\]

 \begin{figure} [ht!]
 \hrule \vspace{0.1cm}
\noindent  {\tt function} ${\tt ADA}\,\big(W,H, \delta, \epsilon_S,\beta\big)$\\
 {\sl computes the new value of  $\beta$.}\\[-0.2cm]
  \hrule
\begin{tabbing}
\hspace*{0.5cm}\=compute $\epsilon_N$ and $\rho$
  according to (\ref{sto3}) and (\ref{sto4});\\[0.1cm]
 \>  {\tt if}\  $\rho<1\  {\tt and}\  \beta>8\  {\tt and}\ (\delta<0.01\  {\tt or}\  \rho<0.8)\  {\tt then}\  \beta
 =\beta/8$,\\[0.1cm]
  \>  {\tt else if}\ $\rho<1\  {\tt and}\  \beta>4\  {\tt and}\ (\delta<0.1\  {\tt or}\  \rho<0.9)\  {\tt then}\  \beta
  =\beta/4$,\\[0.1cm]
  \>  {\tt else if}\ $\rho<1\  {\tt and}\  \beta>2\  {\tt then}\  \beta=\beta/2$,\\[0.1cm]
 \>  {\tt else}\  $\beta =\beta \min\big(8,\rho^2\big);$\\[0.1cm]
  \> {\tt return} $\beta$;\\[-0.7cm]
\end{tabbing}
   \hrule
\caption{\label{f3} The  function ${\tt ADA}$ implements the proposed adaptive strategy. }
\end{figure}

\section{The experimentation}\label{expe}
The experimentation has been performed with a 3.2GHz 8-core Intel Xeon W
processor machine.

\subsection{The datasets}
The experimentation was carried out on both real-world  and artificially generated datasets.
More precisely, we have used for our analysis three classes of matrices.

\medskip
\noindent Class 1: It consists of  matrices of the form $A=VV^T$, where $V\in {\bf R}_+^{n\times p}$ has
random elements uniformly distributed over $[0,1]$.
Namely, three matrices $R1$, $R2$ and $R3$ are generated with $n=2000$ and $p=20,\,40,\,80$.

\medskip
\noindent Class 2: The matrices of this class are obtained starting from undirected  weighted graphs associated
to three real-world datasets of documents ${\bm m}_i$, $i=1,\ldots,n$.
The considered datasets are

\noindent (1)\quad MC: a collection of $n=1033$ medical abstracts from  Medline
(Medical Literature Analysis and Retrieval System Online).

\noindent (2)\quad F575: UMist Faces collection of $n=575$ gray-scale
$112 \times 92$ images  of 20 different people \cite{umist}.

\noindent (3)\quad  F400: Olivetti Faces collection of $n=400$
gray-scale $64 \times 64$ images  of several different people,
     from the Olivetti database at ATT.

\noindent Collection (1) were downloaded  from

http://www.dcs.gla.ac.uk/idom/ir\_resources/test\_collections/

\noindent Collections (2) and (3)  were downloaded  from

http://www.cs.nyu.edu/$\sim$roweis/data.html

\smallskip\noindent
For collection (1) a similarity matrix $A$ is constructed through the  usually considered weights for text data, i.e.
the cosine similarity between two documents
\begin{equation}\label{tsim}
 a_{i,j}=\frac{{\bm m}_i^T{\bm m}_j}{\|{\bm m}_i\|_2
   \|{\bm m}_j\|_2}\,, \quad \hbox{for}\quad i\neq j, \quad \hbox{and}\quad
 a_{i,i}=0.
 \end{equation}
For  collections (2) and (3) a similarity matrix $A$ is constructed using  weights $ e_{i,j}$ suitable for
image data, followed by the normalized cut
\[
 a_{i,j}=d_i^{-1/2}e_{i,j} d_j^{-1/2} \; \hbox{where}\quad
 d_i=\sum_{r=1}^ne_{i,r},\quad \hbox{for}\quad i=1,\ldots,n.
 \]
The weights are expressed through a Gaussian kernel of the form
\begin{equation}\label{isim}
 e_{i,j}=\exp \Big(-\,\frac{\|{\bm m}_i-{\bm
 m}_j\|_2^2}{\sigma^2}\Big)\,,\quad
  \hbox{for}\quad i\neq j,
\quad \hbox{and}\quad e_{i,i}=0,
 \end{equation}
where $\sigma$ is a global scaling
parameter, chosen as the mean value of the distances $\sigma_i$ of the $i$th point ${\bm m}_i$ from its
$7$th nearest neighbor \cite{perona}.

\begin{figure} [ht]
\epsfig{file=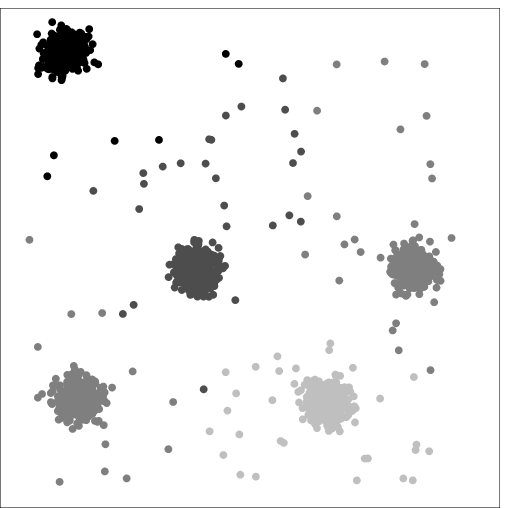,width=2.7cm} \quad
\epsfig{file=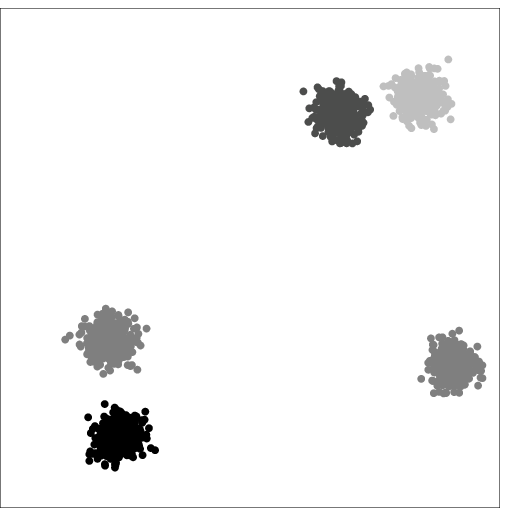,width=2.7cm} \quad
\epsfig{file=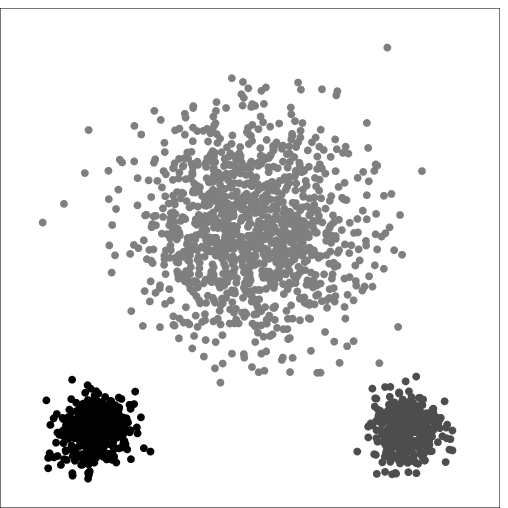,width=2.7cm} \quad
\epsfig{file=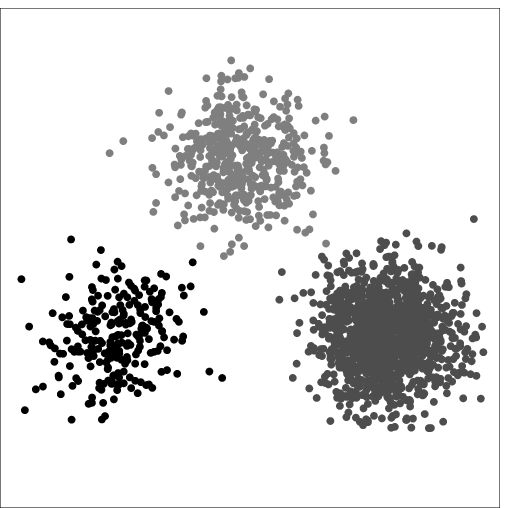,width=2.7cm}
\caption{\label{dat} Synthetic
datasets of points in ${\bf R}^2$: from left to right WSN, SC, SK and DD. }
\end{figure}

\medskip
\noindent Class 3: The matrices of this class are obtained starting from undirected  weighted graphs associated to
four synthetic data sets of points in ${\bf R}^2$ suggested in \cite{liu} (see Figure \ref{dat}):
dataset   WellSeparatedNoise (WSN) consists
of five clusters generated with the same variance and  noise points in the amount of 5$\%$ ; dataset SubClusters (SC)
has three
clusters, and two of them can be divided into subclusters; dataset
SkewDistribution (SK)  has three clusters with different dispersion; dataset  DifferentDensity (DD) has
clusters with different cardinality. All the datasets are generated
with $n=$ 1000, 2000, 4000 and 8000 points, in order to assess the sensitivity of the
algorithm to the increase of the dimensions.
For each dataset the similarity matrix $A$ is constructed as in Class 2 with the choice
 $\sigma=\frac{\sqrt2}{10}\  \max_{i,j} \|{\bm m}_i-{\bm m}_j\|_2 $.

For the matrices of Classes 1 and 2 we look for factors $W$ with  ranks $k=5$, 10, 20, 40 and 80.
For the matrices of Class 3 we look for factors $W$ with  ranks $k=3$, 5 and 10.
By the term ``problem'' we mean a pair $(A_n,k)$ where $A_n$ is the given symmetric matrix of size $n$ and $k$ is the rank of
 the sought matrix $W$.
 Classes 1 and 2 consist of 15 problems each, while Class 3 consists of 48 problems.

\subsection{The tests}
Function {\tt Sym\_ANLS}  calls four update functions:
function {\tt ADA} and, for comparison, three geometrical updatings  with ratio $\zeta=1.01$
(proposed in \cite{kuang1}, here denoted {\tt G1.01}, which gives a slow progression),
$\zeta=1.1$ (here denoted {\tt G1.1}, which gives a mid-level progression) and $\zeta=1.4$ (here denoted {\tt G1.4},
which gives the faster progression).
As {\tt Inner\_solve}, {\tt BPP} and {\tt GCD} are called.
 {\tt GCD} is called with different values of the tolerance $\eta$ used
in the stopping condition, namely
$\eta=10^{-\ell}$, with $\ell=1,\ldots , 5$.
 Due to the fact that in general only approximations of a local minimum of
problem (\ref{sympro}) can be expected, for each problem  and each  instance of function {\tt Sym\_ANLS}, five randomly
generated matrices $W^{(0)}$ have been considered as starting points.
The five runs were performed in parallel using five of the eight available cores
 and the solution with the best final  error has been selected.
The number of outer iterations and the final error of this solution  are indicated as $\nu_{tot}$ and  $\epsilon_S$,
while the largest running time in seconds of the five runs, indicated as $T$, is considered  in order to estimate the
true cost of the whole processing.

\subsubsection{Testing the performance of {\tt ADA}}
The first set of experiments is aimed at evaluating the strategy for updating $\beta$. Table \ref{ta1} shows
the values of  $\epsilon_S$, $\nu_{tot}$ and $T$, averaged on the problems of each class, for the inner methods {\tt BPP}
and {\tt GCD} with $\eta=10^{-3}$.
The averaged results obtained with the other values of $\eta$ are not listed since they,
in comparison with the results of {\tt BPP}, are pretty much the
same of those shown for $\eta=10^{-3}$.

\begin{table} [!h]
\centering
\renewcommand{\arraystretch}{0.3}
\scriptsize
\begin{tabular}{|c|c|ccc|ccc|ccc|}
\hline
&&&&&&&&&&\\
&&\multicolumn{3}{|c|}{Class 1}&\multicolumn{3}{|c|}{Class 2}&\multicolumn{3}{|c|}{Class 3}\\
&&&&&&&&&& \\
\cline{3-11}
&&&&&&&&&& \\
 updating&inner& $\epsilon_S$& $\nu_{tot}$
 &$T$&  $\epsilon_S$
 &$\nu_{tot}$& $T$
 & $\epsilon_S$& $\nu_{tot}$& $T$\\
&&&&&&&&&& \\
\hline
&&&&&&&&&& \\[0.01cm]
&{\tt BPP} &0.010&6.&17.71
&0.416&18.67&8.04& 0.313&13.02&93.01\\[0.01cm]
{\tt ADA}&&&&&&&&&&\\[0.01cm]
&{\tt GCD}&0.010& 16.73&16.96&
0.415& 17.20& 1.62&0.309& 14.90&13.68\\[0.15cm]
\hline
&&&&&&&&&& \\[0.01cm]
&{\tt BPP} &0.010&298.3&775.8
&0.416&18.80&7.36& 0.314&58.94&471.5\\[0.01cm]
{\tt G1.01}&&&&&&&&&&\\[0.01cm]
&{\tt GCD}&0.010& 309.6&195.8&
0.415& 18.40& 1.62&0.309& 60.67& 61.13\\[0.15cm]
\hline
&&&&&&&&&& \\[0.01cm]
&{\tt BPP} &0.010&57.80&164.9
&0.419&17.27&7.10& 0.314&22.94&145.57\\[0.01cm]
{\tt G1.1}&&&&&&&&&&\\[0.01cm]
&{\tt GCD}&0.010& 58.73&37.47&
0.418& 17.27& 1.51&0.309& 23.77& 21.55\\[0.15cm]
\hline
&&&&&&&&&& \\[0.01cm]
&{\tt BPP} &0.010&21.80&68.60
&0.436&11.73&4.16& 0.316&12.27&79.13\\[0.01cm]
{\tt G1.4}&&&&&&&&&&\\[0.01cm]
&{\tt GCD}&0.010& 21.80&14.05&
0.435& 11.87& 1.06&0.309& 12.06& 9.92\\[0.15cm]
\hline
\end{tabular}
\caption{ \label{ta1}  Behavior of ${\tt Sym\_ANLS}$ applied with the inner methods {\tt BPP} and {\tt GCD}
with $\eta=10^{-3}$, averaged on the problems of each class.}
\end{table}
While for each class all the methods appear to be quite equivalent from the point of view of the error, remarkable differences appear from the point of
view of the number of outer iterations and the required time. In general, {\tt BPP} has a smaller number of outer iterations than {\tt GCD}, but
a much larger $T$, indicating that a single outer iteration of {\tt BPP} costs much more than a single outer iteration of {\tt GCD}.
The time comparison shows that, at least in our experimentation, the exact local termination of {\tt BPP} does not pay over
the approximated termination of the iterative method.
For this reason in the following we do not consider {\tt BPP}  anymore.

Turning to the behavior of the geometrical updatings of $\beta^{(\nu)}$, it appears that in general a slower progression
requires more time than a faster progression, with a possible advantage of the error.
As a consequence, it can be very difficult to determine a reasonable ratio of
the updating which combines a low time with an acceptable error level. On the contrary, {\tt ADA} adaptively
produces a dynamical evolution of $\beta^{(\nu)}$
which guarantees on average low computational times and comparable errors.

\begin{figure} [!ht]
\hskip 0.2cm$\epsilon_S^{(\nu)}$\hskip 3.8cm $\epsilon_N^{(\nu)}$\hskip 3.3cm $\beta^{(\nu)}$
\vskip 0.01cm
\epsfig{file=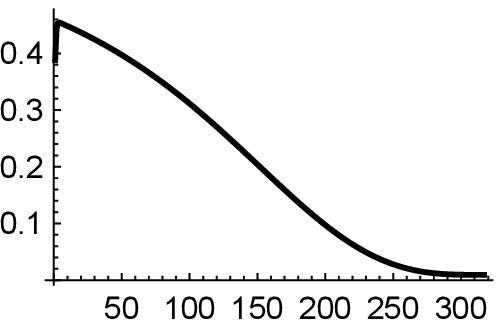,width=3.6cm} \quad
\epsfig{file=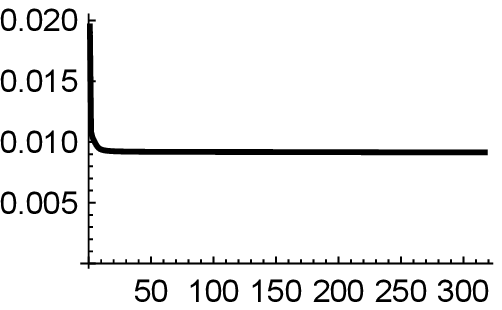,width=3.6cm} \quad
\epsfig{file=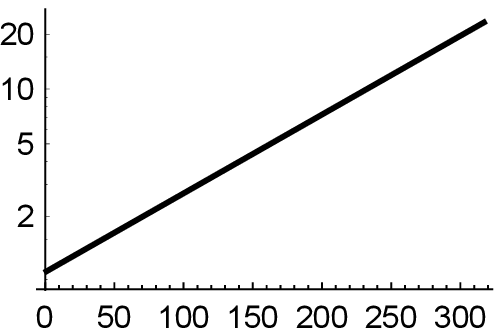,width=3.6cm}
\vskip -0.5cm
\hskip 3.65cm$\nu$\hskip 3.9cm $\nu$\hskip 3.9cm $\nu$
\caption{\label{rot1} Plots of $\epsilon_S^{(\nu)}$ (left), $\epsilon_N^{(\nu)}$ (center) and $\beta^{(\nu)}$ (right)
versus the outer
iteration number for problem R3 with $k=80$. ${\tt Sym\_ANLS}$ is applied with inner method ${\tt GCD}$
with $\eta=10^{-4}$ and function ${\tt ADA}$ (top row) or geometrical updating ${\tt G1.01}$ (bottom row).}
\end{figure}

\begin{figure} [!ht]
\hskip 0.2cm$\epsilon_S^{(\nu)}$\hskip 3.6cm $\epsilon_N^{(\nu)}$\hskip 3.6cm $\beta^{(\nu)}$
\vskip 0.01cm
\epsfig{file=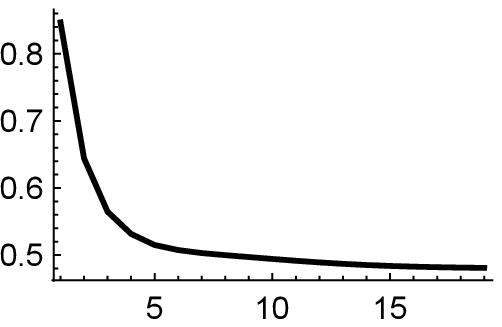,width=3.6cm} \quad
\epsfig{file=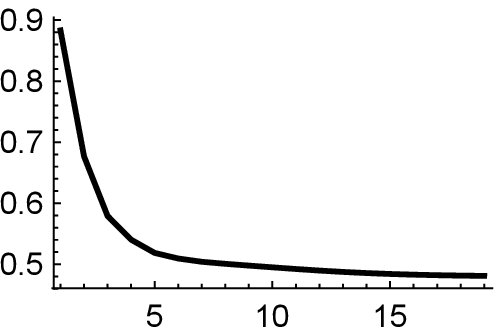,width=3.6cm} \quad
\epsfig{file=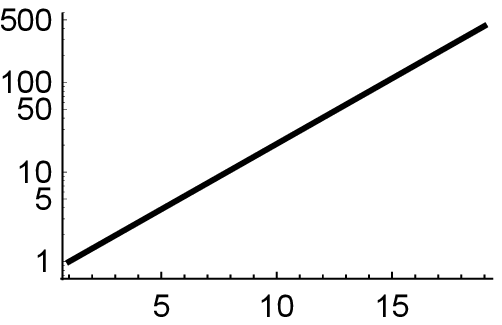,width=3.6cm}
\vskip -0.5cm
\hskip 3.55cm$\nu$\hskip 3.9cm $\nu$\hskip 3.9cm $\nu$

\caption{\label{rot2}  Plots of $\epsilon_S^{(\nu)}$ (left), $\epsilon_N^{(\nu)}$ (center) and $\beta^{(\nu)}$ (right)
versus the
outer iteration number for problem WSN with $n=1000$ and $k=5$. ${\tt Sym\_ANLS}$ is applied with inner method ${\tt GCD}$
with $\eta=10^{-2}$  and function ${\tt ADA}$ (top row) or geometrical updating ${\tt G1.4}$ (bottom row).}
\end{figure}

The  two following examples  present typical situations  where the geometrical updating is outperformed by the adaptive
updating. The first example (see Figure \ref{rot1}) shows how a low rate geometrical updating is outperformed by {\tt ADA}
by the point of view of the cost. The second  example (see Figure \ref{rot2}) shows how  a fast rate geometrical
updating is outperformed by {\tt ADA} by the point of view of the error.

\begin{figure} [!ht]
\hskip 0.2cm$\epsilon_S^{(\nu)}$\hskip 3.7cm $\epsilon_N^{(\nu)}$\hskip 3.6cm $\beta^{(\nu)}$
\vskip 0.02cm
\epsfig{file=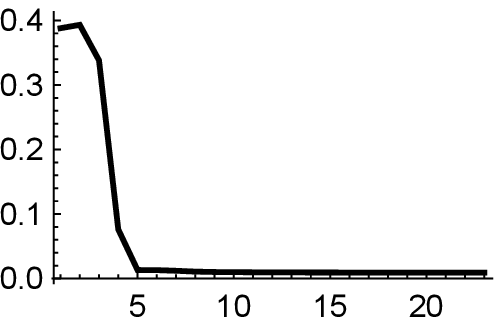,width=3.6cm} \quad
\epsfig{file=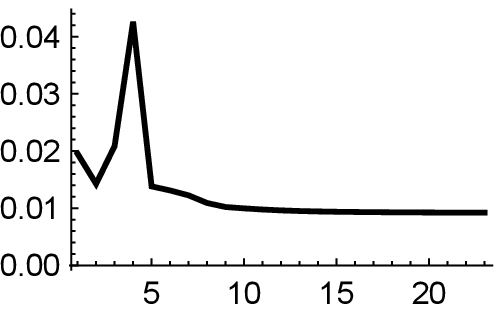,width=3.6cm} \quad
\epsfig{file=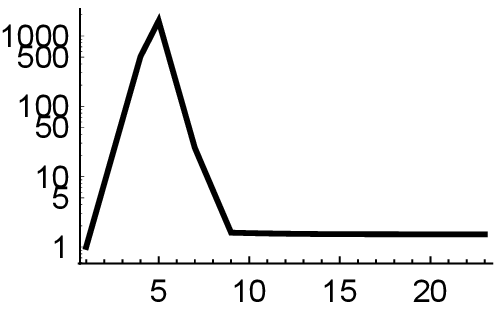,width=3.6cm}
\vskip -0.5cm
\hskip 3.55cm$\nu$\hskip 3.9cm $\nu$\hskip 3.9cm $\nu$
\end{figure}

Figure \ref{rot1}  shows the behaviors of $\epsilon_S^{(\nu)}$, $\epsilon_N^{(\nu)}$  and $\beta^{(\nu)}$
for problem R3 with $k=80$
when    ${\tt Sym\_ANLS}$ is applied  with inner method ${\tt GCD}$ with $\eta=10^{-4}$. The top row plots are
obtained with the updating function ${\tt ADA}$,
the bottom row plots are obtained with the updating function ${\tt G1.01}$.
By inspection of the $\epsilon_S^{(\nu)}$ and $\epsilon_N^{(\nu)}$ plots it appears that  the adaptive updating produces
a transition
from a local minimum of the nonsymmetric error  to another local minimum  of both the nonsymmetric and the symmetric errors.
This transition is obtained through  a fast increase of $\beta^{(\nu)}$  followed by a fast decrease of $\beta^{(\nu)}$
to the value $1.517$.
The final error $\epsilon_S=0.00922$ is obtained in 54.65 sec. with 23 outer iterations.
When  ${\tt Sym\_ANLS}$ is combined with the low rate geometrical updating ${\tt G1.01}$ the
final error $\epsilon_S=0.00921$ is obtained in 655 sec. with 317 outer iterations and a final value $\beta=23.2$,
i.e. a very high cost has to be payed to obtain a comparable error.

Figure \ref{rot2}  shows the behaviors of $\epsilon_S^{(\nu)}$, $\epsilon_N^{(\nu)}$  and $\beta^{(\nu)}$ for  problem WSN
with $n=1000$ and $k=5$,
when  ${\tt Sym\_ANLS}$ is applied  with inner method ${\tt GCD}$ with $\eta=10^{-2}$. The top row plots are obtained with
the updating function ${\tt ADA}$,
the bottom row plots are obtained with the updating function ${\tt G1.4}$.
For this problem  the $\epsilon_S^{(\nu)}$ and $\epsilon_N^{(\nu)}$ plots are very similar from the beginning,
so $\beta^{(\nu)}$
is decreased by ${\tt ADA}$.
It appears that the adaptive updating produces  a transition from a
local  minimum  of the  symmetric error with  $\epsilon_S^{(\nu)}= 0.5$
to another local minimum, assumed as final, with $\epsilon_S=0.28$ in 0.23 sec. with  9 iterations.
This transition is obtained  with a sudden  decrease of $\beta^{(\nu)}$ which reaches and maintains the value $0.5$.
When  ${\tt Sym\_ANLS}$ is combined with the fast rate geometrical updating ${\tt G1.4}$,  $\beta^{(\nu)}$ reaches
the final value  427,
in 0.494 sec. with 19 outer iterations, obtaining the final error  $\epsilon_S= 0.48$ comparable with the one of
the first minimum obtained
by using function ${\tt ADA}$,  i.e.  a worst error is obtained with comparable time.

\begin{figure} [!ht]
\hskip 0.2cm$\epsilon_S^{(\nu)}$\hskip 3.6cm $\epsilon_N^{(\nu)}$\hskip 3.5cm $\beta^{(\nu)}$
\vskip 0.02cm
\epsfig{file=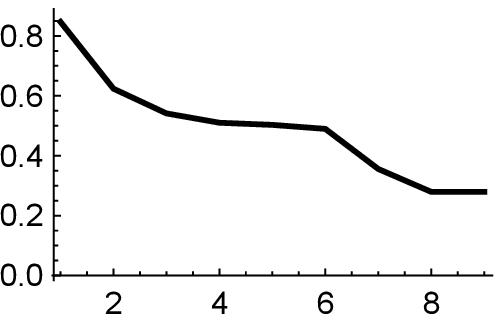,width=3.6cm} \quad
\epsfig{file=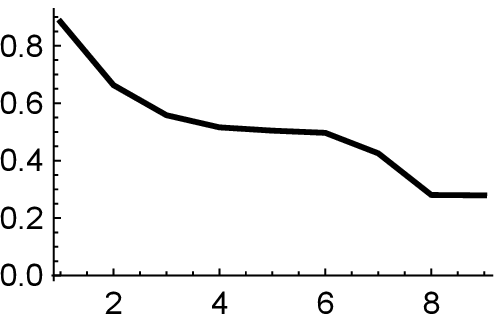,width=3.6cm} \quad
\epsfig{file=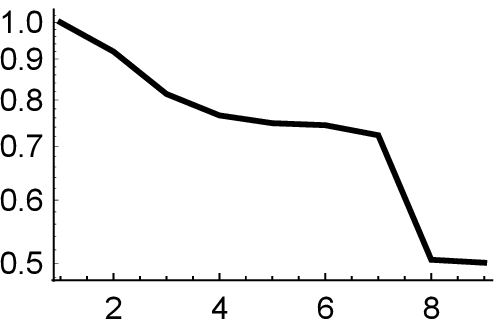,width=3.6cm}
\vskip -0.5cm
\hskip 3.55cm$\nu$\hskip 3.9cm $\nu$\hskip 3.9cm $\nu$

\end{figure}

\subsubsection{Analyzing the performance of {\tt GCD} in the {\tt Sym\_ANLS} schema}
Once {\tt ADA} has been chosen as the most effective updating strategy and {\tt BPP} has been discarded since
more time demanding,
the second set of experiments is aimed at examining how the choice of the tolerance $\eta$ used in the
stopping condition of {\tt GCD} influences the performance. Table \ref{ta2} shows
the values of $\epsilon_S$, $\nu_{tot}$ and $T$, averaged on the problems of each class, for the inner method {\tt GCD}
with $\eta=10^{-\ell}$, $\ell=1,\ldots , 5$.

  \begin{table} [!h]
\centering
\renewcommand{\arraystretch}{0.3}
\scriptsize
\begin{tabular}{|c|ccc|ccc|ccc|}
\hline
&&&&&&&&&\\
&\multicolumn{3}{|c|}{Class 1}&\multicolumn{3}{|c|}{Class 2}&\multicolumn{3}{|c|}{Class 3}\\
&&&&&&&&& \\
\cline{2-10}
&&&&&&&&& \\
 $\eta$&$\epsilon_S$& $\nu_{tot}$
 &$T$& $\epsilon_S$
 &$\nu_{tot}$& $T$
 &$\epsilon_S$& $\nu_{tot}$& $T$\\
&&&&&&&&& \\
\hline
&&&&&&&&& \\[0.01cm]
$10^{-1}$ &0.010&40.8&29.98
&0.419&38.27&4.02& 0.308&18.54&18.49\\[0.15cm]
$10^{-2}$ &0.010&19.60&15.95
&0.416&20.13&1.95& 0.309&14.58&13.71\\[0.15cm]
$10^{-3}$ &0.010&16.73&16.96
&0.415&17.20&1.62& 0.309&14.90&13.68\\[0.15cm]
$10^{-4}$ &0.010&15.33&17.87
&0.416&16.60&1.63& 0.308&14.27&14.11\\[0.15cm]
$10^{-5}$ &0.010&13.80&19.81
&0.415&17.40&1.72& 0.309&14.21&13.64\\[0.15cm]
\hline
\end{tabular}
\caption{ \label{ta2}  Behavior of ${\tt Sym\_ANLS}$ applied with updating {\tt ADA} and inner method {\tt GCD}
with different $\eta$, averaged on the problems of each class.}
\end{table}

In general, a smaller $\eta$ entails a smaller number of outer iterations, but the consequences on the computational time
are not immediate and require a deeper analysis.

At each outer iteration $\nu$ of ${\tt Sym\_ANLS}$ the two matrices $H^{(\nu)}$ and $W^{(\nu)}$ are computed by using {\tt GCD},
which has an initialization phase where the gradient and the Hessian of the objective function and the quantity $\mu$
 of (\ref{mu}) are computed.
After the initialization phase, each inner iteration performs a single coordinate correction.
In a standard implementation, the initialization phase requires a number of  floating point operations ${\mit \Gamma}(k,n)$
 of order $k\, n^2$ and
each coordinate correction requires a number of floating point operations $\gamma(k)$ of order $k$.
Hence for each problem, besides the number of outer iterations $\nu_{tot}$, also the total number $cor$
of single coordinate corrections on both  $H^{(\nu)}$ and $W^{(\nu)}$, $\nu=1,\ldots,\nu_{tot}$, has to be considered.
Of course,
both $\nu_{tot}$ and $cor$ depend on $\eta$, then the overall cost of a run of ${\tt Sym\_ANLS}$ can be expressed as
 \begin {equation}\label{cc}
 c_{tot}(\eta)=c_{out}(\eta)+c_{inn}(\eta),\ \hbox{where}\  c_{out}(\eta)=2\,\nu_{tot}\, {\mit \Gamma}(k,n),
 \  c_{inn}(\eta)= cor\, \gamma(k).
 \end{equation}
In order to analyze how $c_{tot}(\eta)$ depends on the choice of $\eta$ a specific experimentation is made on
three problems, one for each class. The following tables show the behavior of ${\tt Sym\_ANLS}$ applied with updating
{\tt ADA} and inner method {\tt GCD}
with different $\eta$, on the chosen problems. For each problem the solution  with the best final error  has been selected
among the performed five runs.
In  Tables \ref{ta3}, \ref{ta4} and \ref{ta5} the total and inner running times in seconds spent to obtain
this selected solution are denoted  by $T_{tot}$ and $T_{inn}$ (they correspond to $c_{tot}(\eta)$ and $c_{inn}(\eta)$).
The considered problems are R3 of Class 1, MC of Class 2 and WSN with $n=8000$ of Class 3. The values chosen for $k$ are
 $k=80$ for the first two cases and $k=10$ for the third case.
\begin{table} [!h]
\centering
\renewcommand{\arraystretch}{0.3}
\scriptsize
\begin{tabular}{|c|cccccc|}
\hline
&&&&&&\\
 $\eta$&$\epsilon_S$& $\nu_{tot}$ &$cor/1000$&$T_{tot}$&$T_{inn}$
 &$T$\\
&&&&&& \\
\hline
&&&&&& \\[0.01cm]
$10^{-1}$ &0.009&64&8216&105.4&2.46&122.1
\\[0.15cm]
$10^{-2}$ &0.009&24&31666&47.28&9.09&50.75
\\[0.15cm]
$10^{-3}$ &0.009&22&50020&50.20&14.55&51.98
\\[0.15cm]
$10^{-4}$&0.009&22&62609&53.38&18.03&53.38
\\[0.15cm]
$10^{-5}$&0.009&22&91755&61.67&26.28&62.45
\\[0.15cm]
\hline
\end{tabular}
\caption{ \label{ta3}  Behavior of ${\tt Sym\_ANLS}$ applied with updating {\tt ADA} and inner method {\tt GCD}
with different $\eta$, on  problem R3  with $k=80$.}
\end{table}

\begin{table} [!h]
\centering
\renewcommand{\arraystretch}{0.3}
\scriptsize
\begin{tabular}{|c|cccccc|}
\hline
&&&&&&\\
 $\eta$&$\epsilon_S$& $\nu_{tot}$ &$cor/1000$&$T_{tot}$&$T_{inn}$
 &$T$\\
&&&&&& \\
\hline
&&&&&& \\[0.01cm]
$10^{-1}$ &0.360&58&243&24.69&0.221&37.64
\\[0.15cm]
$10^{-2}$ &0.358&27&582&11.55&0.339&11.57
\\[0.15cm]
$10^{-3}$ &0.356&23&1238&10.13&0.595&10.98
\\[0.15cm]
$10^{-4}$&0.356&23&1979&10.38&0.859&11.71
\\[0.15cm]
$10^{-5}$&0.355&23&2856&10.66&1.167&13.42
\\[0.15cm]
\hline
\end{tabular}
\caption{ \label{ta4}  Behavior of ${\tt Sym\_ANLS}$ applied with updating {\tt ADA} and inner method {\tt GCD}
with different $\eta$, on  problem MC  with $k=80$.}
\end{table}

\begin{table} [!h]
\centering
\renewcommand{\arraystretch}{0.3}
\scriptsize
\begin{tabular}{|c|cccccc|}
\hline
&&&&&&\\
 $\eta$&$\epsilon_S$& $\nu_{tot}$ &$cor/1000$&$T_{tot}$&$T_{inn}$
 &$T$\\
&&&&&& \\
\hline
&&&&&& \\[0.01cm]
$10^{-1}$ &0.115&32&132&103.2&0.050&103.2
\\[0.15cm]
$10^{-2}$ &0.115&25&368&82.82&0.051&141.1
\\[0.15cm]
$10^{-3}$ &0.115&19&621&60.45&0.062&73.97
\\[0.15cm]
$10^{-4}$&0.116&13&958&42.14&0.074&54.29
\\[0.15cm]
$10^{-5}$&0.115&17&1520&51.04&0.105&64.26
\\[0.15cm]
\hline
\end{tabular}
\caption{ \label{ta5}  Behavior of ${\tt Sym\_ANLS}$ applied with updating {\tt ADA} and inner method {\tt GCD}
with different $\eta$, on  problem WSN  with  $n=8000$ and $k=10$.}
\end{table}

In the three tables the values of the $T_{inn}$ column are significantly  smaller than  the corresponding values
of the $T_{tot}$ column and, when $k$ is small, the $T_{inn}$ column is negligible compared to the $T_{tot}$ column.
This result
agrees  with the theoretical estimate (\ref{cc}), taking into account the values of $\nu_{tot}$ and $cor$, and
it shows that the inner phase contributes to the cost less than the outer phase.
From Tables \ref{ta3} and \ref{ta4}
it appears that decreasing values of $\eta$
induce a nonincreasing number of outer iterations $\nu_{tot}$ and an increasing number of total corrections $cor$.
As a consequence, also $T_{inn}$ increases. The initial decrease of $\eta$ leads to a decrease  of $T_{tot}$, since the decrease of the outer cost prevails on
the increase of the inner cost. For smaller values of $\eta$, the number of outer iterations
stabilizes leading to an increasing  $T_{tot}$. This behavior is less evident in Table \ref{ta5} which refers to a problem
where a small $k$ is coupled with a much larger $n$.

An analysis of $T_{tot}$ would suggest that  an intermediate value
for $\eta$ appears to be a good choice.
However, the algorithm is called with five different starting points $W^{(0)}$, and the
$T$ column, showing the largest running time cost, represents the effective cost in our
parallel environment. Of course $T$ is greater than $T_{tot}$, but typically shares the same behavior of $T_{tot}$,
and gives the same suggestion for the choice of $\eta$,  confirming what was already shown in Table \ref{ta2}
on average for the problems of each class.

\begin{figure} [!ht]
\hskip 0.2cm$cor_{av}^{(\nu)}$\hskip 3.3cm $cor_{av}^{(\nu)}$\hskip 3.2cm $cor_{av}^{(\nu)}$
\vskip 0.02cm
\epsfig{file=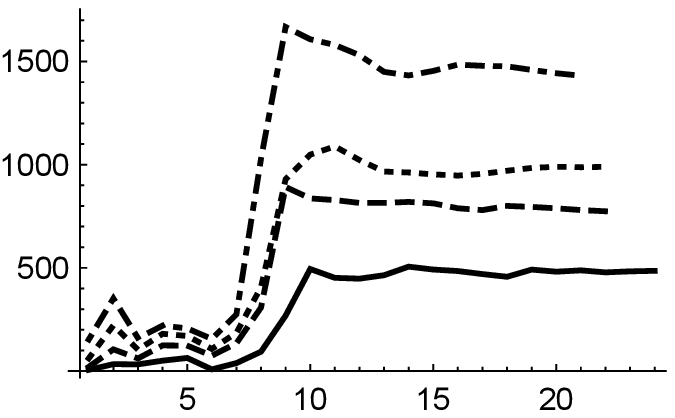,width=3.8cm} \quad
\epsfig{file=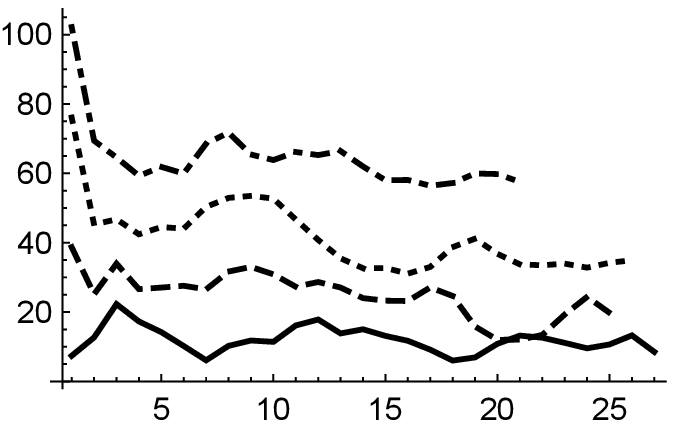,width=3.8cm} \quad
\epsfig{file=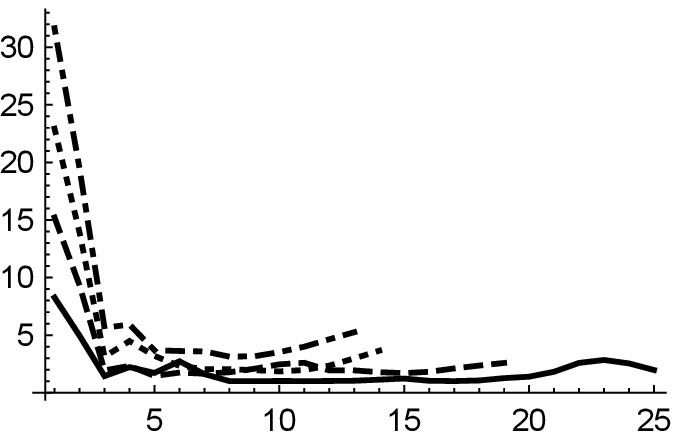,width=3.8cm} \quad
\vskip -0.5cm
\hskip 3.8cm$\nu$\hskip 4.1cm $\nu$\hskip 4.1cm $\nu$
\caption{\label{rot3} Plots of $cor_{av}^{(\nu)}$ for
problem  R3 with $k=80$ (left), problem MC with $k=80$ (middle), problem WSN with $n=8000$ and $k=10$
(right). The inner tolerance is $\eta=10^{-2}$ (continue line), $\eta=10^{-3}$ (dashed line),
$\eta=10^{-4}$ (dotted line), $\eta=10^{-5}$ (dashed-dotted line).}
\end{figure}

A better understanding of the behavior of the inner phase varying $\eta$ can be acquired through Figure \ref{rot3}, where
$cor_{av}^{(\nu)}$ is the number of the coordinate corrections of W and H performed at the $\nu$th outer iteration, divided
by $2n$. For each figure, these average behaviors corresponding to $\eta=10^{-\ell}$ with $\ell=2,\ldots,5$ and
starting with the same initial $W^{(0)}$ are shown.
The figures show in a greater details how the results of the previous tables are formed. In any case we can see that
a larger number of outer iterations corresponds to a lower number of average corrections.
The steep increase of $cor_{av}^{(\nu)}$ in the first figure happens at the same time of the analogous increase of
$\beta^{(\nu)}$ and  corresponds to a change of the local minimum point (see Figure \ref{rot1} top row).

\section{Conclusions}\label{conc}
In this paper an adaptive strategy, called {\tt ADA},  has been introduced  for the updating of the  parameter $\alpha$ in the
penalized nonsymmetric minimization problem (\ref{nonsympro}), when such problem is solved by applying an {\tt ANLS} method.
An extensive experimentation has shown that, when compared to geometrical updatings, {\tt ADA}
produces a dynamical evolution of $\alpha^{(\nu)}$ which guarantees low average computational times and comparable errors.
Moreover, both {\tt BPP} and {\tt GCD} have been tested as inner solvers in the {\tt Sym\_ANLS} schema, concluding
that the latter outperforms the former from the point of view of the computational cost and that
an intermediate value of the internal tolerance  $\eta$, i.e $\eta=10^{-2}$ or $ 10^{-3}$,
should be preferred with {\tt GCD}, especially when $k$ is not very small in comparison with $n$.


\begin{thebibliography}{99}

\bibitem{bjo}
{\AA}. Bj\"{o}rck, {\it Numerical Methods for least squares problems},
SIAM, Philadelphia 1996.

\bibitem{cic}
A. Cichocki and A.H. Phan, {\it Fast local algorithms for large scale nonnegative matrix and tensor factorizations} IEICE Transactions on Fundamentals, E92-A(3) (2009) 708--721.

\bibitem{umist}
D.B. Graham and N.M. Allinson,  {\it Characterizing Virtual
Eigensignatures for General Purpose Face Recognition}, in Face
Recognition: From Theory to Applications ; NATO ASI Series F,
Computer and Systems Sciences, Vol. 163; H. Wechsler, P. J.
Phillips, V. Bruce, F. Fogelman-Soulie and T. S. Huang (eds),
446-456, 1998.

\bibitem{grippo}
L.~Grippo and M.~Sciandrone, {\it On the convergence of the block nonlinear
Gauss-Seidel method under convex constraints}, Oper. Res. Lett., {\bf 26}
(2000) 127--136.

\bibitem{dhil}
C.J. Hsieh and I.S. Dhillon, {\it Fast coordinate descent methods
with variable selection for non-negative matrix factorization},
Proceedings of the 17th ACM SIGKDD International Conference on
Knowledge Discovery and Data Mining  (2011) 1064--1072.

\bibitem{kimpark}
H.~Kim and H.~Park, {\it Nonnegative matrix factorization based on alternating
nonnegativity constrained least squares and active set method}, SIAM J. Matrix
Anal. Appl., {\bf 30} (2008) 713--730.

\bibitem{kimpark2}
H.~Kim and H.~Park, {\it Fast nonnegative matrix factorization: an
active-set-like method and comparisons}, SIAM J. on Scientific Computing, {\bf 33}
(2011) 3261--3281.

\bibitem{kimpark3}
J.~Kim, Y.~He and H.~Park, {\it Algorithms for nonnegative matrix and tensor
factorization: an unified view based on block coordinate descent framework}, J.
Glob Optim, {\bf  58} (2014)  285--319.

\bibitem{kuang1}
D. Kuang D, S. Yun and H. Park, {\it SymNMF: nonnegative low-rank approximation
of a similarity matrix for graph clustering}, J. Glob Optim. {\bf 62} (2015) 545--574.

\bibitem{law}
C.L.~Lawson and R.J.~Hanson, {\it Solving least squares problems},
Prentice-Hall, Englewood Cliffs, N.Y., 1974.

\bibitem{lin}
C.J. Lin, {\it Projected gradient methods for nonnegative matrix factorization}, Neural Computation {\bf 19} (2007) 2756--2779.

\bibitem{liu}
Y. Liu, Z. Li and H. Xiong, {\it Understanding and Enhancement of internal clustering validation measures}, IEEE Trans. Cybernetics, {\bf 43} (2013) 982--993.

\bibitem{paa}
P.~Paatero and U. Tappert, {\it Positive Matrix Factorization: a non-negative
factor model with optimal solution of error estimates of data values},
Environmetrics, {\bf 5} (1994) 111--126.

\bibitem{ppp}
V.P. Pauca, J. Piper and R.J. Plemmons, {\it Nonnegative Matrix Factorization for spectral data analysis}, Linear Algebra and Its Applications, {\bf 416} (2006) 29--47.

\bibitem{perona}
L. Zelnik-Manor  and P. Perona, {\it Self-tuning spectral
clustering}, in Advances in Neural Information Processing Systems,
{\bf 17} (2004) 1601--1608.


\end{thebibliography}
\end{document}